\documentclass[12pt]{amsart}
\usepackage{amscd}
\usepackage[arrow,matrix]{xy}
\usepackage{graphicx, tikz}
\usepackage{hyperref}
\usepackage{comment}
\usepackage{amsmath}

\usepackage{amsmath, latexsym, amssymb}
\numberwithin{equation}{section}
\theoremstyle{plain}
\newtheorem{lemma}{Lemma}[section]
\newtheorem{proposition}[lemma]{Proposition}

\theoremstyle{definition}
\newtheorem{definition}[lemma]{Definition}
\newtheorem*{definition*}{Definition}
\newtheorem{remark}[lemma]{Remark}
\newtheorem{example}[lemma]{Example}
\usepackage{color}
\usepackage{cancel}

\definecolor{brown}{RGB}{150,100,0}

\definecolor{purple}{RGB}{150,0,100}

\definecolor{grey}{RGB}{128,128,128}

\newcommand{\bR}{{\mathbb R}}

\newcommand{\E}{{\mathbb E}}

\newcommand{\N}{{\mathbb N}}

\newcommand{\Bb}{{\mathcal B}}
\newcommand{\Dd}{{\mathcal D}}

\newcommand{\Ff}{{\mathcal F}}
\newcommand{\Mm}{{\mathcal M}}    %moduli space

\newcommand{\Pp}{{\mathcal P}}

\newcommand{\Ss}{{\mathcal S}}

\newcommand{\Xx}{{\mathcal X}}
\newcommand{\Yy}{{\mathcal Y}}

\newcommand{\g}{{\mathfrak g}}

\newcommand{\pb}{{\mathbf p}}

\def\NABLA#1{{\mathop{\nabla\kern-.5ex\lower1ex\hbox{$#1$}}}}
\def\Nabla#1{\nabla\kern-.5ex{}_#1}

\newcommand{\la}{\langle}
\newcommand{\ra}{\rangle}

\renewcommand{\:}{\colon}

\mathchardef\mhyp="2D

\DeclareMathOperator{\Lin}{Lin}

\mathchardef\mhyp="2D

\begin{document}
\title[Natural differentiable structures on statistical models]{Natural differentiable structures on statistical models  and the  Fisher metric}

	\author[H. V. L\^e]{H\^ong V\^an L\^e}	
	\address{Institute  of Mathematics of the Czech Academy of Sciences,
	Zitna 25, 11567  Praha 1, Czech Republic}
\email{hvle@math.cas.cz}

%	\date{\today}

\thanks{Research  of HVL was supported  by  GA\v CR-project 22-00091S and	 RVO:67985840}

\begin{abstract}
In  this paper    I    discuss     the relation   between   the  concept  of the Fisher  metric    and   the concept  of   differentiability  of a    family of       probability  measures.  I 
compare the concepts of smooth statistical manifolds,   differentiable  families of measures,  $k$-integrable parameterized  measure  models, diffeological  statistical  models,    differentiable  measures, which arise  in Information Geometry, mathematical  statistics  and  measure theory, and discuss  some    related problems.	
\end{abstract}

\maketitle
%\tableofcontents
\begin{flushright}
	Dedicated to Professor  Shun-ichi Amari on his 88th  birthday.
\end{flushright}

\section{Introduction}\label{sec1}

In 1945 Rao introduced   the Fisher metric, also called the Fisher-Rao metric,  on    smooth finite dimensional parametric families of probability densities  and,  equipped with this tool, he derived the  Cram\'er-Rao inequality \cite{Rao1945}.  A smooth structure on a   set may give information about its qualitative properties   as well as simplify computations involving the elements of the set.
Between 1970  and 1980   Chentsov and Amari  independently discovered   a one-parameter    family of    affine connections, called {\it $\alpha$-connections}, that are   crucial in   study  of the  geometry  of  smooth finite dimensional parametric families of probability densities, see also  Remark  \ref{rem:chentsov} below.  Using the newly discovered {\it information  geometric structure} consisting of the $\alpha$-connections  and the Fisher metric,  Amari proposed    a differential geometric framework  for constructing  a higher  order  asymptotic  theory  of statistical inference.   Inspired by these  results,   Lauritzen  introduced the concept of (finite dimensional)  statistical manifolds  \footnote{See  Historical Remarks in  \cite[Subsection 1.3]{AJLS17} for more detailed accounts.}.  
In this  article we  re-examine the  relation between  the concept  of the Fisher metric   and the concepts  of smooth parametric  and nonparametric  statistical models, statistical manifolds  and  their generalizations,  which appear  in  measure theory,   functional  analysis,  mathematical statistics  and   Information  Geometry.

The present  article is organized as follows.  In  Section \ref{sec:statm}   we  discuss  the concept
of  a smooth   statistical  model, the  Fisher  metric  and   the  Amari-Chentsov tensor  that has been  proposed  by Amari in his seminal  work \cite{Amari85}, see also   his  book  with  Nagaoka \cite{AN00},   and   independently proposed by Chentsov   \cite{Chentsov72}.  We also recall the related concept of a statistical  manifold, proposed by  Lauritzen \cite{Lauritzen75}.  In Section \ref{sec:difff} we  recall   the concept of  a differentiable   family  of measures, following Bogachev \cite{Bogachev2010},  and relate it to the  concept of  a smooth statistical model,  the concept  of a $k$-integrable parameterized  measure  model, introduced  by Ay-Jost-L\^e-Schwachh\"ofer in
\cite{AJLS15} and \cite{AJLS18},   and the concept of a differentiable  measure   introduced by Fomin \cite{Fomin68}  and by Skorohod \cite{Skorohod1974}.  In  Section \ref{sec:diffl} we   consider    a motivating example  due to Friedrich \cite{Friedrich91},   and   explain  it by  the concept of  a diffeological statistical model  introduced  by L\^e  \cite{Le2020}  that combines   the concept  of a smooth diffeology  due  to Soureau \cite{Souriau1980} and the  concept  of a  differentiable  family  of  measures. We compare  the concept  of  the  Fisher  metric     on $C^k_\tau$-diffeological  statistical models  and  the related  Cr\'amer-Rao inequality   with  existing concepts of  the Fisher  metric  and Cr\'amer-Rao  inequalities in the literature.    
In the      final  section  we  summarize    our discussion  and propose few  questions   for  future research.

\section{Smooth  statistical models  and statistical manifolds} \label{sec:statm}

Motivated   by the local  character  of asymptotic estimation in statistical inference, in 1980
Amari introduced the information geometric structure on  {\it finite dimensional smooth  parametric families    of   probability density functions}  $p (\cdot, \theta)$ on a measurable space	 $ \Xx$  with respect to some  common dominating $\sigma$-finite  nonnegative measure $\mu$  on  $\Xx$ and  parameterized by $\theta$ belonging to some   open  subset  $ \Theta \subset \bR^n$ \cite[p. 2,4,11]{Amari85}, \cite[p.26]{AN00}. The  smoothness of   a  family  of probability density functions $\{p(x, \theta),\, x \in \Xx, \theta \in \Theta\}$   is expressed    via the smoothness   of the function $p (x, \theta)$ in the parameter $\theta\in \Theta$.  Let $\Ff(\Xx)$ denote the set of all real valued functions on $\Xx$.   Amari also required  that the map  $\hat p: \Theta \to \Ff(\Xx), \, \theta \mapsto   p (\cdot, \theta),$ is injective. Amari called  a smooth  family of  probability density  functions  satisfying  the  aforementioned   properties
{\it  an $n$-dimensional statistical model}.  In their influential book  Amari-Nagaoka  also assumed   the following two regularity  conditions (i) and (ii). 
%\begin{itemize}

(i)  The following  rule 	
\begin{equation}\label{eq:dint}
\int _\Xx\partial_V  p(x, \theta)\, d\mu(x)  = \partial_V \int_\Xx p(x, \theta)\,d\mu(x) = 0
\end{equation}
holds for any   tangent  vector $V \in  T_\theta \Theta = \bR^n$  \cite[(2.4)]{AN00}. (This formula  is used for   simplification of  many computations  on   statistical  models.)

(ii)  $p (x, \theta) > 0$  for  all $ x\in \Xx$  and  all $\theta \in \Theta$.
(This  assumption  makes sure that  the set 
${\rm sppt} \,p_\theta := \{ x \in \Xx \vert  p (x, \theta) > 0\} $   does not depend  on $\theta$ \cite[p. 27-28]{AN00},  and hence, all the probability  measures $p(\cdot, \theta)\, \mu \vert \, \theta \in \Theta\}$ are equivalent.
%\end{itemize}

For $V, W \in  T_\theta  \Theta$ we set
\begin{equation}\label{eq:fisher}
\g_\theta(V, W)= \int_{\Xx} \partial _V  \log p (x, \theta) \partial _W  \log p (x, \theta)\, p(x, \theta)d\mu(x).
\end{equation}
Assuming that  the RHS  of \eqref{eq:fisher}    is a finite number  for all $\theta$  and $V, W$,  the   quadratic form    $\g$ is called the {\it Fisher  metric}  on
an $n$-dimensional   statistical  model  $\{ p (\cdot, \theta)\vert\, \theta \in \Theta \subset \bR^n\}$.   Note    we can rewrite 
\eqref{eq:fisher}  as follows
\begin{equation}\label{eq:fisher2}
\g_\theta(V, W)=  4\int_\Xx\partial _V \sqrt{p (x, \theta)}\partial _W \sqrt{p (x, \theta)}\, d\mu(x).
\end{equation}

Formula  \eqref{eq:fisher2} implies  that  the Fisher metric is well-defined  at $\theta$ if   and only if the  function
$\partial _V  \sqrt{p (\cdot, \theta)}$ belongs  to $L^2 (\Xx, \mu)$
for any    $V \in  \bR^n$. 
Now we assume  further   that $\g_\theta$  is positive definite  for all $\theta \in \Theta$.  Under this assumption    the Fisher  metric $\g$ is a Riemannian metric on  $\Theta$. 

On a smooth Riemannian  manifold   $(M, g)$,   one considers the Levi-Civita  connection $\nabla ^{LC} = \nabla ^{LC}(g)$  which  is the unique   torsion-free  metric connection on $M$.   Any other affine  connection   $\nabla$ on $M$  is    differed  from the   Levi-Civita  connection  by   a  tensor   $T^*\in \Gamma(T^*M \otimes T^*M \otimes   TM)$, namely for any $V \in TM$ and  any  vector  field  $X $ on $M$ we have
$$ \nabla   _V   X  = \nabla ^{LC}_V  X +   T^*(V, X) .$$
Clearly, the     connection $\nabla$ is torsion-free  if  and only if $T$ is symmetric.
Amari's {\it $\alpha$-connections}  $\nabla ^\alpha, \alpha \in \bR,$   on   a     statistical model of probability   density functions $\{ p(\cdot,\theta), \, \theta \in \Theta\}$, where the  Fisher metric  $\g$  exists  as a Riemannian  metric, is defined  as follows
\begin{equation*}%\label{eq:amari1}
\nabla ^\alpha =  \nabla ^{LC}(\g)  -\frac{\alpha}{2} T^*.
\end{equation*}
Here    $T^*$   satisfies the following equation  for any $V_1, V_2, V_3 \in T_\theta \Theta = \bR^n$
\begin{equation}\label{eq:AC}
\g (T^* (V_1, V_2), V_3)=\int_\Xx \partial_{V_1}  \log p (x, \theta) \partial_{V_2}  \log p (x, \theta) \partial _{V_3}  \log p (x, \theta)\, p(x, \theta)d\mu(x). 
\end{equation}

(iv) Amari-Nagaoka posed   the  assumption that  the covariant tensor  $T$, defined  by the equation
$T  (V_1, V_2, V_3) =  \g (T^* (V_1, V_2), V_3)$,   is well-defined      for any $\theta\in \Theta$. In this case,  since $T$ is 3-symmetric, $\nabla ^\alpha$  is   a   torsion-free  affine   connection.

\begin{remark}\label{rem:chentsov}
	The   $\alpha$-connections   have been  studied  by  Chentsov    in the case of   discrete  and  finite   sample space  $\Xx$  \cite{Chentsov72}.   Chentsov    proved  that  the Fisher  metric   and   the family of  $\alpha$-connections on  smooth families  of  probability measures  over finite  sample  spaces are  the unique  up to multiplication  Riemannian metric  and  family of connections that are  invariant under  sufficient  statistics \cite[Theorem 11.1, p.159]{Chentsov72}, \cite[Theorems 12.2, 12.3, p. 175, 178]{Chentsov72}.  For a  measurable space $\Xx$  we denote by $\Pp(\Xx)$ the set of all probability measures on $\Xx$.
	Recall    that   a measurable  mapping $f:\Xx \to \Yy$  is  called  a  {\it  sufficient  statistic  for    a  family  $\{\mu_\theta \in \Pp (\Xx) , \theta \in \Theta\}$}, if  there  exists   a Markov kernel  $P: \Yy \times \Sigma _\Xx \to \bR$  such that 
	$P(y \vert \cdot ), y \in \Yy$,  is  the  conditional   probability measure  for  $\mu_\theta$
	relative to $f$  for any $\theta \in \Theta$
	\cite[Definition 2.5, p. 28]{Chentsov72}.  Chentsov's   theorem   has been  generalized  in various forms  by   Ay-Jost-L\^e-Schwachh\"ofer \cite{AJLS15}, \cite{AJLS17}, by L\^e \cite{Le2017}  and  by  Bauer-Bruveris-Michor \cite{BBM2016}.  
	Nowadays   we call the 3-symmetric tensor  $T$ whose values $T(V_1, V_2, V_3)$ is defined in  the RHS  of \eqref{eq:AC}  the 
	{\it Amari-Chentsov  tensor}. 
\end{remark}

Motivated  by the importance  of the Fisher  metric and the  $\alpha$-connections, in  \cite{Lauritzen75}
Lauritzen  proposed  the concept  of   a  {\it statistical manifold} which is a smooth  Riemannian  manifold  $(M, g)$   endowed  with a  3-symmetric covariant tensor  $T$ \cite[Section 4, p. 149]{Lauritzen75}. Lauritzen  asked if  there  is a statistical  manifold  which   does not correspond  to  a  smooth statistical model. In \cite[Section 8.4]{AN00}  Amari-Nagaoka listed   the Lauritzen  question as 
one of mathematical problems posed by  Information Geometry.  Inspired by   Amari-Nagaoka's book and  their list of problems, L\^e gave   an answer  to   Lauritzen's  question in  \cite{Le2005}, see also \cite[Theorem 4.10, p. 222]{AJLS17}.   L\^e's  theorem states that any  smooth statistical manifold $(M, g, T)$ can be immersed  into the  space  $\Pp (\Xx)$  over a countable sample space $\Xx$, which  is  finite if $M$ is compact,  such that  $g$
is induced by the Fisher metric  $\g$ and   $T$ is  induced  by the immersion  from the Amari-Chentsov tensor on   the image  of $M$   in $\Pp(\Xx)$.

In \cite{PS1985} Pistone-Sempi  endowed each  set $\Pp_\mu$  of all probability measures equivalent  to    a given  reference  probability measure $\mu \in \Pp (\Xx)$  with a structure  of an infinite
dimensional smooth  Banach manifold, which can be   provided  with the  Fisher metric  and the  Amari-Chentsov  tensor, see  Example  \ref{ex:smooth}   in  Susection \ref{subs:para}. 
Other infinite  dimensional  statistical  models  that carry a  structure  of a   Banach  or Frech\'et smooth manifold endowed  with the  Fisher  metric have been  considered  recently,   see e.g.  Bauer-Bruveris-Michor \cite{BBM2016},  Newton \cite{Newton2018}, and reference  therein.

\section{Differentiable families of measures,  differentiable measures  and   parameterized measure models}\label{sec:difff}

The  concept  of    a smooth  statistical  model considered  in the previous section  is based  on the  concept  of a   smooth  mapping
from  an open subset   in   a   Fr\'echet space  to   a   space  of  smooth  probability density functions.   In this   section  we shall    consider
a more general   question, what  is  a   differentiable  mapping  from    an open   subset in a topological   vector space  to  the  space  $\Pp (\Xx)$, or  to the space $\Mm(\Xx)$ of all finite nonnegative measures \footnote{We shall     omit  the adjective  ``nonnegative"  to a  measure  in this  paper, and    add ``signed"  to  a measure   if it  is necessary.}    on $\Xx$.  Families   of  finite measures  smoothly  dependent  on a  parameter  arise  in   theory of random  processes, functional analysis, mathematical physics  and mathematical  statistics  see e.g. Bogachev \cite{Bogachev2010}, Borovkov \cite{Borovkov1998}, Chentsov \cite{Chentsov72}, Pfanzagl \cite{Pfanzagl1982}, Pflug \cite{Pflug1996}, Strasser \cite{Strasser1985}.  
All    approaches  to  the concept of  a differentiable  mapping  $f$  into  $\Mm(\Xx)$  exploits  the possibility to approximate  $f$  locally
by    a  linear  mapping into $\Ss(\Xx)\supset \Mm(\Xx)$ and therefore    requires  to specify  a  convergence  type on
$\Ss(\Xx)$.  

\subsection{Convergence types  and natural  topologies     on $\Ss(\Xx)$ and $\Mm (\Xx)$}\label{subs:tau}
%Let $(\Xx, \Sigma_\Xx)$  be  a measurable  space.     Then  
Recall that $\Ss(\Xx)$  endowed with  the total variation norm $TV$  is  a Banach  space  and  the {\it strong   topology}  $\tau_v$   generated by this norm is  compatible with the linear  structure on $\Ss(\Xx)$. It is well-known that  the  {\it weak  topology} $\tau_W$   on the Banach space $\Ss(\Xx)_{TV}$   is also  compatible   with the linear structure  on $\Ss(\Xx)$.  Besides  these topologies, one considers  also   {\it setwise convergence} and  the associated  topology $\tau_s$ which is the weakest topology  on $\Ss(\Xx)$ such that for any $A \in \Sigma_\Xx$  the map  $I_A: \Ss(\Xx) \to \bR, \mu \mapsto  \mu (A),$ is continuous. Equivalently, $\tau_s$  is generated  by  the duality   with  the space $\Ff^s (\Xx)$  of all  simple  functions on $\Xx$. Note  that    a sequence  of  $\mu_n \in \Ss (\Xx)$  converges   setwise  if and only if   $\mu_n$ converge in the weak  topology $\tau_W$, see e.g. Bogachev \cite[Corollary 4.7.26]{Bogachev2010}.
Furthermore, the  restriction    of $\tau_s$ to $\Pp(\Xx)$ is not metrizable  unless   $\Xx$ is countable, see e.g. Ghosal-van der Vaart \cite[p. 513]{GV17}.

If $\Xx$ is  a topological   space  then we   consider   the  Borel  $\sigma$-algebra $\Bb (\Xx)$ unless  otherwise stated.  Some time  we also consider  the Baire $\sigma$-algebra $\Bb a (\Xx)$ that is the smallest  sub-$\sigma$-algebra of $\Bb(\Xx)$ such that  any continuous function on $\Xx$ is measurable.    If $\Xx$ is a  metrizable topological
space,  then  $\Bb (\Xx) = \Bb a (\Xx)$ by Bogachev \cite[Corollary 6.3.5, vol. 2,
p. 13]{Bogachev2007}. The space  $C_b  (\Xx)$  consisting  of bounded  continuous    functions   on a topological  space $\Xx$  is   a  Banach    space with  the sup-norm  $\|f \|_\infty : = \sup_{x \in \Xx} \vert f(x)\vert$.   Clearly  the canonical  embedding $C_b (\Xx)_{\infty}\to L^\infty (\Xx, \mu)$  is continuous for any $\mu \in \Mm (\Xx)$.  Hence $C_b (\Xx)_{\infty}$  is  a  closed  subspace  of    the Banach space  $\Ss (\Xx)'$.  
%\begin{definition}\label{def:weak2}  
We recall that the   {\it weak*  topology}  $\tau_w$  on the space $\Ss(\Xx)$,  where $\Xx$ is a topological space endowed  with the Borel $\sigma$-algebra, %$\Mm_{\sigma} (\Xx)$} of    signed aire  measures 
is  the  weakest topology such that   for each  $ f \in   C_b (\Xx)$   the  map
$I_f : \Ss (\Xx) \to \bR, \:  \mu \mapsto  \int _\Xx f d\mu$, is  continuous.
We also denote by $\tau_w$ the restriction of  $\tau_w$ on  subsets of $\Ss(\Xx)$. %The  weak  topology
%on the space is  the  topology  induced  from $\tau_w$ on $\Mm_{\sigma} (\Xx)$.
%\end{definition}
It is well-known that if $X$ is infinite, then $(S(X), \tau_w)$ is non-metrizable  but the restriction  of $\tau_w$ to $\Mm (\Xx)$ is metriazable, if $\Xx$ is a separable metrizable space, see e.g. Bogachev \cite[p. 102]{Bogachev2018}.

\subsection{Differentiable families   of measures  and differentiable  measures}\label{subs:fam}
Let us  begin with the  concept   of  a  diffferentiable family  of signed finite measures, following Bogachev \cite[Section 11.2]{Bogachev2010}, which encompasses  all  notions  of differentiable  families  of measures
that  have  been  considered  in   literature in mathematical statistics   to  the best of our  knowledge,  see   e.g. Chentsov \cite{Chentsov72}, Pfanzagl \cite{Pfanzagl1982},  Strasser \cite{Strasser1985},   Pflug \cite{Pflug1996},  van der Vaart \cite{Vaart1991},  and references  therein. Let $\tau$  be a  topology  on a    vector space  $\Ss(\Xx)$.   

\begin{definition}\label{def:taud}(\cite[Defintions 11.2.1, 11.2.11]{Bogachev2010})
	(i) A family  of signed  finite measures   $\{\mu(t)\in \Ss(\Xx)\vert \, t \in  (a, b)\}$ is called  {\it $\tau$-differentiable}   at a  point $\tau_0 \in (a, b)$  if a limit
	$$\mu' (t_0) = \lim _{ s \to 0} \frac{\mu (t_0 +s) -\mu(t_0)}{s}$$
	exists  in the topology $\tau$.
	
	(ii)   The definition of   continuity is completely  analogous.  Higher order  differentiability
	is defined  inductively. \footnote{See  \cite[Definitions  4.1.10,  4.1.11, p. 428]{BS2017}  and   Example  \ref{ex:ckdi}  below.}
	
	(iii) Let $p\in [1, +\infty)$.   A family of finite nonnegative measures  $\{\mu(t)\in \Mm(\Xx)\vert \,t \in  (a, b)\}$, is called  $L^p(\nu)$-differentiable, or  just {\it $L^p $-differentiable}, at a point $t_0$ if there  exists  a nonnegative $\sigma$-finite measure $\nu$ such that  $\mu (t) = f(t)\, \nu$  for all  $t$ and  the mapping  
	$\Phi:  I \to  L^p (\nu), \, t \mapsto f(t) ^{1/p},$ is differentiable at $t_0$. A family of  signed finite measures  is called $L^p$-differentiable, if  its positive  and negative parts in the  Jordan-Hahn decomposition have this properties.
\end{definition}

%By  Remark \ref{rem:lp}  below, the concept  of $L^p(\nu)$-differentiability of  a finite  dimensional family  of nonnegative  measures     does not depend  on the choice   of $\nu$  and  therefore  we  just write $L^p$-differentiability, omitting $\nu$.

\begin{remark}\label{rem:geom}    
	Similarly, we define  the differentiability  of  a family  of measures $\{\mu(m)\vert \, m \in M\}$ where $M$ is a   (possibly infinite dimensional)  differentiable  manifold, provided     that some concept
	of differentiability of mappings from $M$ to  $\Ss(\Xx)$  is chosen (e.g.  if $M$ is a normed space, or $M$ is  an open subset  of a   topological  vector space,   and the  differentiability can be  G\^ateaux-differentiability   or Hadamard  differentiability,  see also  Bogachev-Smolyanov \cite[Section 4.1]{BS2017}  for  the most general  concept  of  differentiability with  respect  to a system of sets.)
\end{remark}
\begin{remark}\label{rem:amari1}  If a  family   $\{\mu_t\in \Pp (\Xx)\vert  \, t \in (a, b)\}$  is $\tau_s$-differentiable, then  for any $t_0 \in (a, b)$ we have
	$$\int_\Xx d \mu'(t_0)  = 0.$$
\end{remark}

We assume in the remainder  of  this Subsection that  in the case  $\Xx$ is a  topological  space  then $\Xx$ is a  completely regular  topological  space.  Under this assumption,     we have \cite[p. 379]{Bogachev2010}
\begin{equation}\| \mu\| = \sup_{f \in C_b(\Xx),\,  \| f\|_{\infty} \le  1} \int _\Xx  f(x) d\mu(x).\label{eq:sup}
\end{equation}
Equality \eqref{eq:sup}  allows us to   formulate   strong results concerning   $\tau_w$-differentiability  and compare  it with  $\tau_s$- and $\tau_v$-differentiability %  below. %where $\|  f \| : = \sup _x  \vert f(x)\vert$.
in   Proposition \ref{prop:logder} below, which  is    a combination of  \cite[Lemma 11.2.3, Theorems 11.2.5,  11.2.6]{Bogachev2010} and some  assertions in their  proofs. Note  that  \eqref{eq:sup} also holds if we  replace $C_b(\Xx)$ by $\Ff^s(\Xx)$. Furthermore, the $\tau_W$-differentibility is weaker  than the $\tau_v$-differentiability  and    stronger  than the  $\tau_s$-differentibility   and         shall not be considered  in   the present  Subsection.

\begin{proposition}\label{prop:logder}  Let $-\infty <  a < b < \infty$. 
	(1) A family $\{ \mu(t)\in \Mm (	\Xx),\, t\in (a, b)\}$  is $\tau_w$-differentiable   at $t_0 \in (a, b)$  if and only if  it  is   $\tau_s$-differentiable   at $t_0$  and if and only if
	for any  $A \in \Bb (\Xx)$ the function  $t \mapsto \mu_1 (A)$ is differentiable  at $t_0$.
	
	(2) Assume that a family  $\{ \mu(t)\in \Mm (	\Xx),\, t\in (a, b)\}$  is $\tau_w$-differentiable.
	Then the family  $\{\mu_t\}$ is $\tau_v$-continuous, $\mu'(t) \ll \mu(t)$  for all $t \in (a, b)$, moreover  the family  $\{\mu_t\}$ is $\tau_v$-differentiable  for almost  every $t \in (a, b)$.  Furthermore there  exists  $\nu \in \Pp (\Xx)$ such  that  $\mu(t)\ll \nu$ for all $t$. Suppose also  that   the function  $t \mapsto  \|\mu'\ (t)\|$  is Lebesgue integrable on $(a, b)$.  Then  for all $t \in (a, b)$ we have
	\begin{equation}\label{eq:lebesgue}
	\mu(t) -\mu(a) = \int_a ^t  \mu' (t)\, dt
	\end{equation}
	where the integral   on the right  is an $\Ss(\Xx)$-valued  Bochner integral.
	If the $\{ \mu'(t)\}$ is $\tau_v$-continuous, e.g. 
	if $\{ \mu(t)\}$ is twice $\tau_s$- or  twice  $\tau_w$-differentiable, then $\{\mu(t)\}$ is $\tau_v$-differentiable  and  the corresponding    derivative coincides   with the $\tau_s$-derivative or $\tau_w$-derivative   respectively.
\end{proposition}

\begin{remark}\label{rem:dominate}  %(1)
	The proof of  $\mu'(t) \ll \mu(t)$  is  given in the  proof of \cite[Theorem 11.2.6]{Bogachev2010}, see  also \cite[p.142]{AJLS17},  \cite[Lemma 2]{LT2021}.   
	The  Radon-Nykod\'ym  derivative  $\frac{d\mu'(t)}{d\mu(t)}\in L^1 (\mu (t))$  is called  the {\it logarithmic  derivative of $\mu (t)$  at  $t$}  \cite[p. 383]{Bogachev2010}, \cite[(4.3)]{AJLS17} \cite[(3.66), (3.68)]{AJLS18}.
\end{remark}

\begin{remark}\label{rem:domi2} Proposition \ref{prop:logder}  implies  that
	a  $\tau_s$ or $\tau_w$-differentiable  family  $\{ \mu(t)\in \Mm(\Xx), t \in (a, b)\}$  is  dominated   by    a probability  measure  $\nu \in \Pp(\Xx)$,   so  we can  write  $\mu (t) =  f(t, \cdot)\, \nu$  where $f(t, \cdot) \in L^1 (\nu)$.  Now assume that    the function  $t \mapsto  \| \mu'(t)\|$ is Lebesgue  integrable. By the Lebesgue  differentiation theorem, \eqref{eq:lebesgue}  implies that   the  partial  derivative   $\partial _t  f (t, \cdot)$ exists  for  a.e. $t \in (a, b)$  and
	\begin{equation}\label{eq:lebesgue0}
	\partial _t  f(t, \cdot) =  \frac{d\mu'(t)}{d\mu(t)}\in L^1 (\mu (t)) %\partial _t \log  \mu (t)  
	\end{equation}
	as an equality  for  elements in $L^1 (\mu (t))$.
	
In \cite[Definition 3.5]{AJLS17}  a  family of  dominated  measures $\{ f(\theta, \cdot)\, \mu\vert\, \theta \in \Theta\}$ is called {\it regular}, if   for all $\theta \in \Theta$ the partial   derivative $\partial _\theta  f(\theta, \cdot)$   exists  and belongs to $ L^1 (\mu)$.
\end{remark}

\begin{example}\label{ex:amari}  Finite  dimensional statistical models $\{p (\cdot, \theta)\vert, \theta \in \Theta \subset \bR^n \}$    in the previous  section       are    regular   $\tau_v$-differentiable families  of   dominated  measures $\{\theta \mapsto  p (\cdot, \theta)\mu, \theta \in \Theta\}$, where,  by  Proposition \ref{prop:logder}, we can assume that  $\mu \in \Pp (\Xx)$. The  condition  that    $\partial _\theta  f(\theta, \cdot)$  belongs to $ L^1 (\mu)$      is  necessary   for  the existence  of  the LHS of \eqref{eq:dint}.  The  equality \eqref{eq:dint}   follows then  from Remark \ref{rem:amari1}. Furthermore the Fisher  metric  is well-defined   on    a    statistical model  $\{p (\cdot, \theta)\vert, \theta \in \Theta \subset \bR^n \}$    if and only  if  the family  is $L^2$-differentiable.
\end{example}

\begin{proposition}\cite[Proposition 11.2.12 and  its  proof]{Bogachev2010}\label{prop:lp}
	(i)	If   for  some $p \ge 1$   a family $\{ \mu(t)\in \Mm (\Xx)\vert\, t \in (a, b)\}$ is  $L^p(\nu)$-differentiable  at $t_0$, then it is   $\tau_v$-differentiable  at   this point   and its logarithmic derivative  $\frac{d \mu' (t_0)}{d \mu (t_0)}$ exists  and belongs  to   $L^p (\mu(t_0))$. % In this case
	%\begin{equation}\label{eq:logder}
	%\Phi' (t_0) =  p^{-1}f(t_0) ^{1/p} \partial _t \log \mu (t).
	%\end{equation}
	(ii)	Conversely, if the family $\{ \mu(t)\in \Mm(\Xx)\vert\, t \in (a, b)\}$    is  $\tau_s$-differentiable or $\tau_w$-differentiable  and the  function  $\| \frac{d \mu' (t)}{d \mu (t)}\|_{L^p (\mu(t))}$ is locally integrable on $(a, b)$,
	then  this family  is $L^p(\nu)$-integrable  for  some  $\nu \in \Pp (\Xx)$ and  for any $ (a_1, b_1) \subset (a, b)$ we have
	\begin{equation}\label{eq:lebesgue1}
	f(b_1) ^{1/p} - f(a_1)^{1/p} = p ^{-1} \int_a ^b f(t) ^{1/p} \frac{d\mu'(t)}{d \mu (t)} \mu(t)\,  dt
	\end{equation}
	$\nu$-a.e. and	as an equality for $L^p(\nu)$-valued Bochner integral. 
\end{proposition}
%We    also  write  \eqref{eq:lebesgue1}   in the following  equivalent  way, cf. \eqref{eq:lebesgue}
%\begin{equation}\label{eq:lebesgue1a}
%f(b_1) ^{1/p}\nu - f(a_1)^{1/p} \nu= p ^{-1} \int_a ^b f(t) ^{1/p}  \mu'(t)\,  dt
%\end{equation}
%as a  Bochner integral.

\begin{remark}\label{rem:lp}   Proposition  \ref{prop:lp}(i)  implies that  a family $\{ \mu(t)\in \Mm (\Xx)\, t \in (a, b))\}$  is $L^p(\nu)$-differentiable, if and only if  it is  $L^p (\nu')$-differentiable      where  $\nu' \in \Pp (\Xx)$  dominates   $\mu(t)$ for all  $t\in (a, b)$. The same argument  implies  that, given  a smooth  finite dimensional manifold  $M$,  a family  $\{\mu(y)\in \Mm (\Xx), \,  y \in M\}$  is  $L^p(\nu)$-differentiable, if and only if  it is  $L^p (\nu')$-differentiable      where  $\nu' \in \Pp (\Xx)$  dominates   $\mu(y)$ for all  $y \in M$.  % Therefore  one uses the notion  of the  $L^p$-differentiability, omitting $\nu$.
\end{remark}

\begin{example}\label{ex:diffmeasure} Given a measurable mapping $\kappa: \Xx \to \Yy$  between measurable spaces $\Xx , \Yy$,  we denote  by $\kappa_*: \Ss(\Xx) \to \Ss(\Yy)$ the pushforward   mapping  defined by $\kappa_* \mu   (B): =\mu (\kappa^{-1} (B))$ for $\mu \in \Ss(\Xx)$ and $B \in \Sigma_\Yy$.  Now let $\Xx$ be a   linear  space endowed  with a  $\sigma$-algebra $\Sigma_\Xx$  that  is invariant under  a shift $L_{th}: \Xx \to \Xx, x \mapsto x - th$  for   some $h \in \Xx$ and $t \in \bR$.  In particular, the shift  $L_{th}: \Xx \to \Xx$ is  a measurable  mapping.  Then   a finite  signed  measure  $\mu \in \Ss(X)$  is called {\it  differentiable along  the vector $h$ in Fomin's sense}, if  the family  $(L_{th})_* \mu$  is  $\tau_s$-differentiable \cite{Fomin68}. 
Bogachev   \cite[p. 69-70]{Bogachev2010} proposed an equivalent definition  of   Fomin's   differentiability  of $\mu\in \Ss(\Xx)$  that requires  a  weaker condition, namely   for  every  $A \in \Sigma_\Xx$  there  exists  a finite  limit
	$$d_h\mu (A): = \lim_{t \to 0} \frac{\mu (A +th) -\mu(A)}{t}.$$
A measure  on $\Bb a(\Xx)$ is called  a  Baire measure.  A Baire measure  $\mu$ on a  locally convex  topological  vector  space $\Xx$ is called  {\it Skorohod  differentiable  along  a vector $h \in \Xx$},  if for   every  $f \in C_b(\Xx)$  the function
	$$ t \mapsto  \int_\Xx f (x -th) d \mu (x)$$
	is  differentiable \cite{Skorohod1974}, \cite[Definition 3.1.5, p.71]{Bogachev2010}. Skorohod  differentiability  of a Borel measure  is understood  as the   differentiability  of its restriction to $\Bb a(\Xx)$. Note that  every     topological  vector space  is  completely regular \cite[Theorem 1.6.5, p. 44]{BS2017}. Using  a  theorem due to Alexandroff one can show that    $\mu$ is  Skorohod   differentiable along $h$, if and only if    the family $ (L_{th})_* \mu$  is  $\tau_w$-differentiable \cite[p.71]{Bogachev2010}. \footnote{The concept  of $\tau_w$-differentiability  in Definition \ref{def:taud}  extends  natural to   families  of Baire  measures.}
\end{example}

\subsection{Parameterized measure models}\label{subs:para}
In \cite{AJLS17}, \cite{AJLS18}   Ay-Jost-L\^e-Schwachh\"ofer  proposed a
geometric  formulation  of  the $L^p$-differentiability  of a   family  of finite signed measures   by   using the concept  of the $p$-th root  of a finite nonnegative measure.   Assume that $ 1\le p \in \bR$. For $\mu \in \Mm (\Xx)$ we   set
$$\Ss^{1/p} (\Xx, \mu): = \{  \nu \in \Ss (\Xx, \mu)\vert \, \frac{ d\nu}{d\mu} \in L^p ( \mu)\},$$
$$\Ss^{1/p}_0 (\Xx, \mu) = \{ \nu \in \Ss^{1/p} (\Xx, \mu)\vert \, \int_\Xx d v = 0\}.$$
The natural identification $\Ss^{1/p} (\Xx, \mu) = L^p ( \mu)$ defines   a  $p$-norm  on  $\Ss^{1/p} (\Xx, \mu)$ by setting
$$ \| f \mu\| _{ p} : = \|  f \| _{ L^ p  (\mu)}.$$
Then  $\Ss^{1/p} (\Xx, \mu) _{p}$  is a Banach space.
For   $\mu_1 \ll \mu_2$     the linear inclusion 
$$\Ss ^{1/p}  (\Xx, \mu_1)\to \Ss^{1/p}  (\Xx, \mu_2), \:   f \mu_1  \mapsto  f (\frac{d \mu_1}{d\mu_2})^{1/p} \in \Ss^{1/p} (\Xx, \mu_2)$$
preserves the $p$-norm. 
Since $(\Mm (\Xx), \ll)$ is a directed  set,  the  directed  limit  
\begin{equation}\label{eq:limr}
\Ss^{1/p} (\Xx)_{p} : = \lim _{\longrightarrow} \{\Ss^{1/p} (\Xx,\mu)_{p}\vert \; \mu \in \Mm (\Xx)\}
\end{equation}
is a Banach space. % Clearly   the zero  measure  $\{0\}$ belongs  to $\Ss ^{1/p} (\Xx, \mu)$  for  any
%$\mu \in \Mm  (\Xx)$. If  $\mu \in \Mm ^* (\Xx) : = \Mm (\Xx)\setminus \{0\}$ then   
%$$\mu^{1/p}  : =\frac{\mu}{\|\mu\| ^{p}} \in \Ss^{1/p} (\Xx, \frac{\mu}{\|\mu\|})_{p} $$
%for any  $p \ge 1$. 
The  image  of $\mu  \in \Ss^{1/p}(\Xx, \mu)_{p} $ in $\Ss^{1/p}  (\Xx)_{p}$ via the directed limit in \eqref{eq:limr}  is called  {\it the  $p$-th root  of $\mu$}  and denoted   by $\mu ^{1/p}$. %\footnote{This  definition of the $\mu^{1/p}$  is  equivalent  to the definition  of $\mu^{1/p}$ via  the  isomorphism $\Ss ^{1/p} (\Xx, \mu) =  L^p (\Xx, \mu)$ and  taking the   direct limit  of $L^p (\Xx, \mu)$ \cite[p.144]{AJLS17}.}
%The $p$-th root of  $\mu \in \Mm (\Xx)$  is defined  as follows 
%$$\mu ^{1/p}: =  \|\mu\|_{TV}  ^{1/p}\big ( \frac{\mu}{\|\mu\|_{TV}}\big ) ^{1/p}.$$
By \cite[Proposition 3.2]{AJLS17}  the map  $\pi^{1/p}: \Mm (\Xx) \to \Ss^{1/p} (\Xx)_p,  \mu \mapsto \mu^{1/p},$ is continuous  with respect  to  the strong  topology $\tau_v$ on  $\Mm(\Xx)$. % and on the Banach space $\Ss^{1/p} (\Xx)$.
Using the   Jordan-Hahn decomposition  $\mu = \mu^+ -\mu^-$, we  can  extend  the continuous map $\pi^{1/p}$ to   continuous  maps  $\pi^{1/p}_{\pm}: \Ss(\Xx)_{TV} \to \Ss^{1/p}(\Xx)_p, \, \mu \mapsto  (\mu^+)  ^{1/p}\pm (\mu^-) ^{1/p}$ \cite[p. 148]{AJLS17}.   Note that $\pi^{1/p}_- (\Ss (\Xx)_{TV}) = \Ss^{1/p} (\Xx)_p$.    Thus  the  image  $\pi_-^{1/p}(\mu)$  for $\mu \in \Ss (\Xx)$  is   called    the  $p$-th root of measure  $\mu$. This  agrees  with  the concept of
$L^p$-differentiability  of  a family of  signed finite  measures \cite[Definition 11.2.11. p.386]{Bogachev2010}.

\

$\bullet$ Given    topological  vector  spaces  $V, W$  we denote  by  $\Lin (V, W)$  the space  of all continuous  linear maps  from  $V$  to  $W$. If  $V$  and $W$  are Banach spaces   then  $\Lin (V, W)$ is the  Banach space  of   bounded  linear  operators from $V$ to $W$, which is    endowed  with the operator  norm. 

\begin{definition} \cite[Definitions 3.4,  3.7]{AJLS17} (i) A   {\it parameterized  measure  model} 
	is a triple  $(M, \Xx, \pb)$, where  $M$ is a Banach manifold  and  $\pb: M \to  \Mm(\Xx) \stackrel{i}{\to}  \Ss(\Xx)_{TV}$ is a  $C^1$-map, i.e.  the composition $i\circ \pb: M \to \Ss(\Xx)_{TV}$ is  continuously Fr\'echet differentiable. 
	(ii)	 Let $k \ge 1$.  A parameterized measure  model $(M, \Xx, \pb)$  is called    {\it $k$-integrable} 
	if the map $\pb ^{1/k} : = \pi^{1/k} \circ \pb: M \to  S^{1/k} (\Xx)_{k}$  is   continuously  Fr\'echet differentiable.   It is  called
	{\it weakly $k$-integrable}, if  $\pi^{1/k} \circ \pb$  is
	{\it weakly    Fr\'echet differentiable}   and its  derivative is  {\it weakly continuous}.  In other  words,  on any    coordinate chart $(U_j, \varphi_j: U_j \to  E_j)$, where  $U_j \subset M$ and $E_j$ is a Banach space  modeling $M$, for  any  $z \in \varphi_j(U_j)$ there exists  a  {\it weak  differential}
	$d _z(i \circ \pb^{1/k} \circ \varphi_j ^{-1})_ \in \Lin (E_j, \Ss^{1/k} (\Xx)_k)$   such  that  
	$$ \frac{i \circ \pb^{1/k} \circ \varphi _j^{-1} (z + v)  - i \circ \pb^{1/k} \circ \varphi_j ^{-1} (z) - d _z(i \circ \pb^{1/k} \circ \varphi_j ^{-1}) (v)  }{\|v\|_{E_j}} $$
	converges  to $0\in \Ss^{1/k}  (\Xx)_k$  in  the  weak topology  as  $v$ converges  to $0\in E_j$,  moreover,   the map  $(\varphi_j (U_j)) \to \Lin (E_j, \Ss^{1/k} (\Xx)_k), \,    z \mapsto d_z(i \circ \pb^{1/k} \circ \varphi_j ^{-1}),$
	is    weakly continuous.  %$\tau_W$-differentiable and its  derivative is $\tau_W$-continuous.  Furthermore   we call  the model  (weakly) $\infty$-integrable if it  is  (weakly)  $k$-integrable  for all $k \ge  1$.
	
	(iii)  A parameterized  measure  model $(M, \Xx, \pb)$  is called {\it a  parameterized statistical model}, if  $\pb (M) \subset \Pp (\Xx)$.
	
	(iv) A parameterized  measure  model $(M, \Xx, \pb)$  is called  {\it dominated}, if   there exists
	$\mu \in \Pp (\Xx)$ such that  $\pb (m)\ll \mu$ for all $m \in M$.
\end{definition}

By   Proposition \ref{prop:logder}, see  also  Remarks \ref{rem:domi2} and \ref{rem:lp}, we  obtain  the following  Lemma  immediately.

\begin{lemma}\label{lem:equi}   Let $M$ be a finite   dimensional smooth manifold and $k \ge 1$.
	Then  any  $k$-integrable   measure model  $(M,\Xx, \pb)$
	corresponds  to a  continuously $L^k$-differentiable  family  with parameters in $M$  and any  continuously $L^k$-differentiable   family with parameters in $M$
	corresponds  to a $k$-integrable  measure  model $(M, \Xx, \pb)$.	
	%(i) Any  continuously   $L^p(\nu)$-differentiable  family $\{p (y)\nu\in \Mm(\Xx),\,  y \in  M, p: M \to  L^p (\nu)\}$   corresponds  to   a   $p$-integrable  parameterized	measure  model $(M, \Xx, \pb)$ where     $\pb= p \,\nu: M \to \Ss^{1/p} (\Xx, \nu) \subset  S^{1/p}  (\Xx) $.

\end{lemma}

The    value  of the  (weak) differential $d_m \pb^{1/k} \in  \Lin (T_m M, \Ss^{1/k} (\Xx)_k )$  at $V\in T_mM$ is called     the  (weak)   derivative   of $\pb^{1/k}$  along  $V$. For an   arbitrary (weakly) $k$-integrable    parameterized measure  model $(M, \Xx, \pb)$   the  (weak)   derivative of $ \pb^{1/k}$  along  $V$  is given  as   \cite[Proposition 3.4, p. 154]{AJLS17}, cf. Proposition \ref{prop:lp} (ii) 
\begin{equation}\label{eq:derp}
\partial _V \pb ^{1/k} (m) = \frac{1}{k}\frac{ d\partial _V \pb (m)}{d \pb (m)}\pb^{1/k} (m) \text{  for } m \in M, V \in T_m M.
\end{equation} 
The  following  Proposition  is a version of  Proposition \ref{prop:lp}(ii).  
\begin{proposition}\label{prop:kint}\cite[Theorem 3.2, p. 155]{AJLS17}  Let $(M, \Xx, \pb)$  be a parameterized  measure   model.  The  model is $k$-integrable, if and only if  the map $V \mapsto  \|\partial _V   \pb^{1/k}\|_p$ is well-defined  on $TM$   and continuous. The model
	is weakly  $k$-integrable, if   for any $V\in  TM$    the  weak derivative    $\partial_V  \pb^{1/k}$  exists and weakly converges  to $\partial_{V_0}\pb^{1/k}$  as
	$V$ converges  to $V_0$.
\end{proposition}

By Proposition \ref{prop:kint}, any  2-integrable parameterized  statistical model  $(M, \Xx, \pb)$  is   endowed  with the  continuous  Fisher  metric  $\g $  on $M$ defined  at $m \in M$ by
\begin{equation}\label{eq:fisherpara}
\g_m (V, W) = \la  \frac{d\partial _V  \pb (m)}{d\pb (m)} , \frac{d \partial _W \pb (m)}{d \pb (m)}\ra_{L^2 (\pb (m))}.
\end{equation}

\begin{example}\label{ex:smooth}  Assume  that  a  finite    dimensional   statistical  model of  probability  density functions  $\{p (\cdot, \theta) \mu\vert\, \theta \in \Theta, \mu \in \Pp (\Xx)\}$  in    the previous   section  is endowed  with  continuous  Fisher  metric   in  \eqref{eq:fisher2}. Then  it is   a  $L^2$-differentiable family  of probability  measures (Example \ref{ex:amari}),  and  by  Lemma \ref{lem:equi}  it is a
2-integrable  parameterized  statistical model.  
	Furthermore,  the Fisher  metric $\g$  in \eqref{eq:fisherpara}     can be  expressed  as in  \eqref{eq:fisher}.
	The  continuous Amari-Chentsov  tensor  on    a  3-integrable  parameterized meaure model   is defined in a similar way  \cite[(3.96), p. 167]{AJLS17},  cf. \eqref{eq:AC}. 
	
	The Pistone-Sempi  manifold $\Pp_\mu$  is an $\infty$-integrable  parameterized  statistical  model  $(\Pp_\mu, \Xx, i)$, where  $i: \Pp_\mu \to \Pp (\Xx)$ is  the natural inclusion  \cite[Proposition 3.14]{AJLS17}.
	
	Newtons'   statistical manifolds  of differentiable densities  in \cite{Newton2018}
	are  3-integrable    parameterized  statistical  models.
\end{example}

\section{Diffeological  statistical  models}\label{sec:diffl}
The  concept  of  smooth parametric   families  of probability  density  functions considered by Amari as  well as    its  various generalizations         we considered  in the previous  sections  is  based  on  the notion of a smooth map  from a nice    smooth space  $M$  to the  vector    space  $\Ss(\Xx)$. This  concept  leads  to the notion  of  a  smooth   parameterization  of  a  subset  $S\subset \Ss(\Xx)$ by  $M$.    Diffeology   theory  is  created  by Souriau in 1980 \cite{Souriau1980}      to describe  {\it  consistent} smooth  parameterizations  of  a  set. 
Diffeology language is   therefore  a  convenient language  to deal with various   objects  in parametric and nonparametric   Information Geometry
in a  unified   elegant  framework. 

\subsection{Friedrich's example}

In \cite{Friedrich91}, motivated  by Amari's  work \cite{Amari85} and  the concept  of differentiable curves in mathematical statistics \cite{Strasser1985}, Friedrich  supplied   the space 
$P(\lambda):=\{  \mu \in \Pp(\Xx)\vert \, \mu \ll   \lambda\}$, where $\lambda$  is a  $\sigma$-finite  nonnegative measure on $\Xx$, with  the  Fisher metric   as follows.  He    defined   the tangent  space  $T_\mu P (\lambda)$   to be  the  linear space $\Ss_0 ^{1/2} (\Xx, \mu)$  and  the Fisher
metric  on $\Ss^{1/2}_0 (\Xx, \mu)$  is    given  by  (cf. \eqref{eq:fisherpara})
\begin{equation}\label{eq:ffisher}
\g_\mu (V, W) =\big \la \frac{d V}{d \mu},  \frac{d W}{d\mu}\big\ra _{ L^2 (\mu)}.
\end{equation}

Note that  the tangent spaces  $\Ss^{1/2}_0 (\mu)$  and  $\Ss^{1/2}_0 (\mu')$   are not isomorphic  if  $\mu$ and $\mu'$ are not equivalent. Hence $P(\lambda)$  does not have  the structure  of an infinite  dimensional
Fr\'echet  manifold.  There is no obvious  way to  parameterize   $P(\lambda)$  by a single map
from an  infinite dimensional  smooth  manifold. 
In  \cite{Le2020} L\^e proposed   to  endow   $P(\lambda)$   with a natural   diffeology,  called  statistical diffeology, which    provides $P(\lambda)$ with   a       differentiable structure, whose  tangent  space    at $\mu \in P (\lambda)$ is  $\Ss^{1/2}_0 (\mu)$.

\subsection{Diffeological  spaces}

Let us  recall the concept  of  diffeology, following  Iglesias-Zemmour \cite[1.5, 1.14]{IZ13}.

\begin{definition}\label{def:diffeo}
	Let $X$ be a nonempty set.  A {\it  parametrization}  of $X$ is a map from   an open subset $U \subset \bR^n$ to  $X$.  A {\it diffeology} of $X$  is any set $\Dd$ of parametrizations  of $X$
	such that  the following axioms are satisfied
	
	D1.  {\it Covering}  The set  $\Dd$ contains the constant  parametrization  ${\bf x}:  r \mapsto   x$, defined on $\bf R^n$  for all $x \in X$ and   all $n \in \N$.
	
	D2. {\it Locality}   Let $\pb : U \to X$ be a parameterization. If for every point $u \in U$ there exists  an open  neighborhood  $V\subset U$ of $u$ such that $\pb_{\vert V} \in \Dd$, then the parameterization $\pb$ belongs to $\Dd$.
	
	D3. {\it Smooth compatibility}.  For every  element $\pb : U \to \Xx$ of $\Dd$, for every  open subset  $V \subset \bf R^m$  and for every $F\in C^\infty (U, V)$, the composition    $\pb \circ F$ belongs to $\Dd$.

	A {\it diffeological space}   is  a nonempty set  $X$ equipped with a diffeology  $\Dd$.
	A map  $f: X \to X'$  between  two diffeological   spaces $(X, \Dd)$ and  $(X', \Dd')$  is said to be {\it smooth},  if for any  $\pb \in \Dd$ we have  $f \circ \pb \in \Dd'$.
\end{definition}

We verify immediately that the composition  of  smooth  maps between   diffeological spaces
is a smooth map. So the set of all  diffeological spaces     forms a category whose  morphisms are
smooth mappings.

\begin{example}\label{ex:ckdi}        Let $\tau$    be a topology on $\Ss(\Xx)$  defined  in Subsection \ref{subs:tau}  and $U$  an open subset in $\bR^n$.   A      map $f: U  \to (\Ss (\Xx), \tau)$  will be called   {\it $\tau$-differentiable  at $x \in U$},  if  there  exists
	a {\it $\tau$-differential}  $d_x f \in \Lin (\bR ^n,  (\Ss (\Xx), \tau))$ such that
	$$ \frac{f (x +v) - f(x) - d_x f  (v)}{\|  v\|} \stackrel{ \tau}{\to}  0 \in \Ss (\Xx),$$
	as  $v$  goes to  $ 0 \in \bR^n$.  Here $\| v\|$ is the  Euclidean norm of $v$,  and $\stackrel{\tau}{\to} 0$  denotes the  convergence  to  $0 \in \Ss (\Xx)$ in the $\tau$-topology.  If  $f$ is   $\tau$-differentiable  at  all  $ x\in U$, it will be called  {\it $\tau$-differentiable}.
	A  $\tau$-differentiable  map  $f: U \to \Ss (\Xx)$ will be   called   {\it  continuously  $\tau$-differentiable}, or $C^1_\tau$-differentiable,  if   its  {\it partial $\tau$-derivative map}
	$\partial f: TU \to  (\Ss (\Xx), \tau), (x, v) \mapsto  \partial  _v  f  : = d_x f (v)\in  \Ss (\Xx)$ is     continuous  in $\tau$-topology. Inductively   we define the concept  of   $C^k_\tau$-map    from $U$  to $\Ss (\Xx)$.
	Note  that  for  any   Fr\'echet differentiable  map $f: U \to E$, where $E$  is a  Banach  space, the   differential $df: U \to  \Lin (\bR^n, E)$  is continuous if and only if  its  partial derivative map $\partial  f: TU \to  E, (m, v)\mapsto    d_m f (v),$  is continuous. Denote by $ \Lin ^2 (\bR^n, \bR^n; E)$  the set  of  all continuous  bilinear mappings from  $\bR^n \times \bR^n$ to  $E$.  Taking  into account  the equality
	$$\Lin  (\bR^n, \Lin  (\bR^n, E))  =  \Lin ^2 (\bR^n, \bR^n; E)$$
	see  e.g. \cite[Proposition 2.4, p. 7]{Lang99}, we conclude  that $C^k_{\tau_v}$-differentiable maps  from $U$ to  $\Ss (\Xx)$    are  exactly    Fr\'echet  $C^k$-differentiable  maps.
	
	Now    we     denote  by  $\Dd^k_\tau$     the   set  consisting   of  all  $C^k_\tau$ maps  from  $U \subset \bR^n$  to  $\Ss (\Xx)$, where  $U$    runs  over all open  subsets  in $\bR^n$, $ n \in \N$. It is  straightforward to check that   $\Dd^k_\tau$  is a   diffeology   on $\Ss(\Xx)$.
\end{example}

\begin{example}\label{ex:ckdii}    Let $X$ be a set  and  $(X', \Dd')$ be a diffeological   space.  Given a  map $f: X \to X'$   there  is a  coarset diffeology     of $X$ such that  the map $f$ is smooth. This diffeology is called  the {\it  pullback  diffeology}  of the diffeology $\Dd'$ by $f$, and denoted  by   $f^*(\Dd')$   \cite[Section 1.26, p. 14]{IZ13}.  We have  
	$$f^* (\Dd') : = \{ \pb : U \to X\vert \, f \circ  \pb \in \Dd'\}.$$
	If $X$ is a subset of  $X'$  and $f$ is the natural inclusion the  pullback diffeology on $X$ is  also called the subset diffeology \cite[Section 1.33, p. 18]{IZ13}.
\end{example}

\begin{example}\label{ex:ckdiii}  Let $(X, \Dd)$  be a diffeological  space  and $f: X \to X'$  a map.  The  finest diffeology  of $X'$   such that $f$  is smooth  is called   the {\it pushforward  diffeology} of $\Dd$ by $f$, and denoted by $f_* (\Dd)$  \cite[Section 1.43, p. 24]{IZ13}.  A mapping $\pb : U \to X'$  belongs  to $f_* (\Dd)$  if and only if  it satisfies the  following  condition.  For every point $u \in U$  there exists  an open neighborhood $V \subset U$ of $u$  such that   either 
	$\pb _{\vert V}$ is a constant  map, or  there exists  a map ${\bf q}: V \to X$ in  $\Dd$ such that 
	$\pb_{\vert V} =  f \circ \bf q$. 
\end{example}

\begin{remark}\label{rem:diff}     In \cite{Le2020}  L\^e  introduced  the   concept  of    a $C^k$-diffeology,  which    coincides  with the concept of      diffeology in Definition  \ref{def:diffeo}, except  that the    condition     D3  on smooth compatibility  in Definition \ref{def:diffeo} is replaced  by  the $C^k$-compatibility,   namely    the   set $C^\infty(U, V)$ is replaced by the larger set  $C^k (U, V)$.  
	Since
	$C^\infty  (U,V)$ is a  subset of $C^k(U, V)$ any  $C^k$-diffeological  space   is a diffeological space.	
\end{remark}
\subsection{Diffeological   statistical models  and diffeological Fisher metric}

From now  on,    a {\it statistical model}  is a  subset $\Pp_\Xx \subset \Pp (\Xx)$ \cite{McCullagh02}.  Let $\tau$ be a  topology on $\Ss(\Xx)$, with    specification as  in    Subsection \ref{subs:fam}.

\begin{definition}\label{def:diffstat}  Let $\Pp_\Xx\subset \Pp(\Xx)$  be a statistical  model  and $i: \Pp_\Xx \to \Ss(\Xx)$  the inclusion.  A diffeology   $\Dd$ on $\Pp_\Xx$  is   called a {\it  $C^k_{\tau}$  statistical diffeology},   if  the   embedding $i: (\Pp_\Xx, \Dd)\to  (\Ss (\Xx),  \Dd^k_{\tau})$ is a smooth  embedding. In this case   $(\Pp_\Xx, \Dd)$  is called a {\it $C^k_{\tau}$ diffeological statistical model}.
\end{definition}

\begin{example}\label{ex:i}   Any  statistical model  $\Pp_\Xx\stackrel{i}{\to}\Ss(\Xx)$      has  the subset diffeology $i^* (\Dd^k_{\tau})$,  which is  a  $C^k_\tau$ statistical  diffeology.
\end{example}	

\begin{example} \label{ex:iii}  Let  us consider all possible  open subsets  $U \subset \bR^n$, $n \in \N$  and  for  $ p \ge 1$  we  set
	$$\Dd^{1, p}_{\tau_v}: =\{ f: U \to \Pp(\Xx), \vert    \pi^{1/p}\circ  f: U \to \Ss^{1/p}(\Xx)_p  \text{  is  a  Frech\'et $C^1$-map }\}.$$
	Equivalently, by Lemma \ref{lem:equi}, $f$ is a continuously $L^p$-differentiable  map.      Then we verify immediately that the set $\Dd^{1, p}_{\tau_v}$ is a $C^1_{\tau_v}$  statistical diffeology of  $\Pp(\Xx)$. %Furthermore, the  pullback $i ^* (\Dd^2_{can})$ to  $P(\lambda)$ via the canonical inclusion   $P (\lambda) \stackrel{i}{\to} \Mm (\Xx)$  is a subset of $\Dd_2$.
\end{example}

\begin{example}\label{ex:v}  For  any  smooth  Banach  manifold  $M$ we denote  by $\Dd^\infty_{M}$ the    diffeology consisting    of  all  smooth  maps  from      all  open subsets  $U \subset \bR^n$, $n \in \N$,   to  $M$. Let $(M, \Xx, \pb)$ be a   $k$-integrable    parametrized measure  model.   Then    the image  $\pb (M)\subset \Mm(\Xx)$ has   the pushforward diffeology  $\pb_* (\Dd^\infty_{M})$  which is a  $C^1_{\tau_v}$ statistical  diffeology. %The  $k$-integrability of the parameterized measure model implies that  $\pb^* (\Dd_{can(M)})\subset i^* (\Dd^k_{can})$.   
\end{example}

\begin{definition}\label{def:tangent} \cite[Definition 1]{Le2020}
	(i) Let  $(\Pp_\Xx, \Dd)$ be a $C^k_\tau$ diffeological  statistical model.	A   vector $v\in \Ss(\Xx)$   is  called    a {\it tangent  vector }  of   $(\Pp_\Xx, \Dd)$  at $\mu \in \Pp_\Xx$  if  there  is a    smooth map $ c : (-1, 1) \to  (\Pp_\Xx, \Dd)$ such that   $ c(0) =  \mu$ and  $ c' (0) = v$.

	(ii) {\em The tangent  cone}  $C_\mu (\Pp (\Xx), \Dd)$   consists  of all   tangent  vectors $v$ of
	$(\Pp (\Xx),\Dd)$ at   $\mu$.
	
	(iii)  {\em The tangent  space}  $T_\mu  (\Pp (\Xx), \Dd)$  is the linear hull of  $C_\mu (\Pp (\Xx), \Dd)$.
	
\end{definition}

\begin{example}\label{ex:tangentD1}
	Let  us  revisit  Friedrich's  example  of  the space $P (\lambda)$	 where $\lambda$ is a   $\sigma$-finite
	nonnegative  measure on  $\Xx$.  Denote  by  $i : P (\lambda) \to \Ss (\Xx)$   the natural  inclusion.
	By Remark \ref{rem:amari1}  and Proposition \ref{prop:logder},  the tangent space  $T_\mu (P (\lambda), i^* (\Dd^1_{\tau_v}))$ is a subset of  $\Ss^1_0 (\mu)$.
	We shall  show that $T_\mu (P (\lambda), \Dd^1_{can}) = \Ss^1_0(\mu)$.
	Let $\phi\mu \in \Ss^1_0 (\mu)$.
	
	Set $\tilde \mu (t) = p (x, t)\, \mu$, where
	\begin{equation}\label{eq:tangentp1}
	p(x, t) = \begin{cases}   1 +   t\phi (x) &  \text{ if  }   t \phi (x) \ge 0\\
	\exp (t \phi (x)) & \text{  if }  t \phi (x) < 0.
	\end{cases}
	\end{equation}
	In \cite[p. 142]{AJLS17}  Ay-Jost-L\^e-Schwachh\"ofer showed that   the curve
	\begin{equation}\label{eq:tangentpb}
	\mu (t) = \frac{\tilde \mu(t)}{\|\tilde \mu(t)\|}
	\end{equation}
	belongs  to $\Ss^1_0 (\mu)$  is continuously  $\tau_v$-differentiable, and   $\mu(0) = \mu$, $\mu'(0) = \phi \mu$.
	Thus  $T_\mu (P (\lambda), i^* ( \Dd^1_{\tau_v})) = \Ss^1_0 (\mu)$.
	
	Now  let   us consider  the difeology  $i^* (\Dd  ^{1, 2}_{\tau_v})$ of $\Pp (\lambda)$. 
	By Lemma \ref{lem:equi}, Remark \ref{rem:amari1}, and Proposition  \ref{prop:lp},     any  tangent  vector  in  $T_\mu (P (\lambda), i^* (\Dd^{1,2}_\tau))$ belongs  to  $\Ss^{1/2}_0  (\Xx, \mu)$.  We  claim   that     $T_\mu  (P (\lambda ), i^*  (\Dd ^{1, 2}_{\tau _v})) = \Ss^{1/2}_0 (\Xx, \mu)$.   Let  $ \phi\mu \in \Ss ^{1/2}_0  (\Xx, \mu)$.   We  consider  the   same     curve   $\tilde  \mu  (t) = p   (x, t) \mu$, where  $p (x, t)$ is defined  by \eqref{eq:tangentp1}  and  set  $\mu (t)$  by \eqref{eq:tangentpb}, cf. \cite[p. 147]{AJLS17}.
	By Proposition \ref{prop:kint}, $\mu(t)$  is continuously $L^2$-differentiable,   and $ \partial_ t \mu ^{1/2}(0)= {1\over 2}\phi  \mu ^{1/2} $  by \eqref{eq:derp}.
	This   proves  our claim.
\end{example}

\begin{definition}\label{def:dfisher} \cite[Definition 4]{Le2020}.  A  $C^k_\tau$ diffeological statistical model $(\Pp_\Xx, \Dd)$  is called {\it almost 2-integrable},
	if   for  all   $\mu \in \Pp_\Xx$ we have  
	$T_\mu  (\Pp_\Xx, \Dd)\subset \Ss^{1/2}_0 (\mu)$.
	A $C^k_\tau$-diffeological  statistical model  $(\Pp_\Xx, \Dd)$ is called  {\it  2-integrable}, if    the inclusion  $i: (\Pp_\Xx, \Dd)\to (\Pp(\Xx), \Dd^{1,2}_{\tau_v})$ is smooth. 
\end{definition}

Clearly   any   2-integrable $C^k_\tau$-statistical  model  is almost 2-integrable but there is an almost 2-integrable  statistical model which is not  2-integrable, see \cite[Example 4]{LT2021}.

Now  we define   {\it the diffeological   Fisher  metric}  on    an  almost 2-integrable    $C^k_\tau$-diffeological statistical  model  $(\Pp_\Xx, \Dd)$ by the   same formula in \eqref{eq:ffisher}.  Namely  for  $V, W \in  T_\mu (\Pp_\Xx, \Dd)\subset \Ss^{1/2}_0 (\Xx, \mu)$ (the inclusion  is a  consequence  of Proposition \ref{prop:lp} and Remark \ref{rem:amari1})  we  set
$$ \g_\mu (V, W) = \la \frac{dV}{d \mu}, \frac{dW}{d\mu} \ra _{L^2 (\mu)}.$$

\begin{example}\label{ex:2integrd}  Let  $(M, \Xx, \pb)$ be a 2-integrable   parameterized statistical model.
	Then   $(\pb (M), \pb_*  (\Dd^\infty_{M}))$  is a 2-integrable  $C^1_{\tau_v}$-statistical model.    The tangent
	space   $T_{\pb  (m)} (\pb  (M), \pb_* (\Dd^\infty_{M}))$ consists  of  tangent  vectors $ \partial _V\pb  (m), V \in T_m  M$. Then  we have
	$$\g_{\pb (m) }  (\partial  _V \pb (m), \partial_W  (\pb (m)) ) =  \la  \frac{d\partial _V  \pb (m)}{d\pb (m)} , \frac{d \partial _W \pb (m)}{d \pb (m)}\ra_{L^2 (\pb (m))}.   $$
	Thus  the Fisher  metric $\g_m (V, W)$ defined on the  parameter space  $M$ by  \eqref{eq:fisherpara} is the pull back  of  the  diffeological  Fisher  metric in \eqref{eq:ffisher} on  the   image $(\pb (M), \pb_*  (\Dd^\infty_{M}))$  and its  degeneracy  is caused  by  the  kernel  of the  differential $d_m \pb: TM \to \Ss^{1/2}_0 (\Xx, \pb (m))_2$.
\end{example} 

\begin{remark}\label{rem:lt}
	Since  the   diffeological  Fisher  metric $\g$  is   nondegenerate,  L\^e-Tuzhilin     used    $\g$  to define the Fisher  distance  on  2-integrable  $C^k_\tau$-diffeological  statistical  models  $(\Pp_\Xx, \Dd)$.  The length  of a  smooth curve 
	$c: [0,1] \to (\Pp_\Xx, \Dd)$ is   defined  by \cite[Definition 5]{LT2021}
	$$ L(c)  = \int_0 ^1  \vert    c ' (t)\vert _\g, dt $$ 
	The length of a piece-wise  smooth curve  is set to be
	the sum  of the lengths of  its smooth sub-intervals.  The diffeological  Fisher distance  $d_\g$ between two-points   is defined  as the infimum of the lengths  over the space of  piece-wise \cite{LT2021}.   L\^e-Tuzhilin   showed   that   a   2-integrable $C^k_\tau$ diffeological    statistical model  $(\Pp_\Xx, \Dd)$ endowed  with the Fisher  distance $d_\g$ is a length space \cite[Theorem 1]{LT2021}, moreover $d_\g (x, y) \ge \| x -y\|_{TV}$ \cite[Lemma 4]{LT2021}.  Thus   the topology  on $\Pp_\Xx$ generated  by   $d_\g$  is  not weaker than   the  strong  topology  $\tau_v$.  If the Hausdorff dimension of $(\Pp_\Xx, \Dd, d_\g)$  is finite, then  $(\Pp_\Xx, \Dd, d_\g)$ is endowed  with  the  Hausdorff-Jeffrey  measure, which coincides  with the  unnormalized   Jeffrey  prior measure defined on   smooth  $n$-dimensional statistical  models  \cite[Theorem  3]{LT2021}.
\end{remark}  

\subsection{Diffeological and parametric  Cram\'er-Rao inequalities}

In density estimation problems, given  a statistical  model  $\Pp_\Xx \subset  \Pp(\Xx)$, we  wish  to     measure   the   accuracy of  a nonparametric {\it estimator  $\hat \sigma: \Xx \to \Pp_\Xx$}, or  its  $\varphi$-coordinate, formalized as  a map $\varphi: \Pp_\Xx \to V$, where  $V$ is a topological vector space.  
If $\Pp_\Xx  = \pb (M)$ where
$(M, \Xx, \pb)$ is a   parameterized   statistical model,   then  it is convenient  to have   a {\it  parametric  estimator} $\hat \sigma: \Xx \to M$, or  its  
$\varphi$-coordinate,  given a map $\varphi: M \to V$.  
%$$\MSE^\varphi _\mu[\hat \sigma] (l, l): = \int_\Xx\big (l \circ \varphi\circ \hat \sigma (x)  - l \circ \varphi \circ \mu )^2\, d\mu (x),$$
%\end{equation*}
%$$\MSE_m  ^\varphi[\hat \sigma](l, l) = \int_\Xx  (l \circ \varphi\circ \hat \sigma (x) - l \circ \varphi \circ  m)^2\, d {\pb (m)} (x), $$
%for nonparametric  and  parametric $\varphi$-estimators, respectively. Here  we assume  that  $\MSE^\varphi _\mu[\hat \sigma] (l, l)$ and $\MSE_m  ^\varphi[\hat \sigma](l, l)$  are well-defined  for all  $\mu \in \Pp_\Xx$ 
%and  $m \in M$, respectively.
We  measure  the   deviation  of a  parametric  and  nonparametric $\varphi$-estimator  $\varphi \circ \hat \sigma: \Xx \to V$  from its  mean  value using  its covariance
$$V_\mu ^\varphi  [\hat \sigma] (l, l) : = \int_\Xx\big (l \circ \varphi\circ \hat \sigma (x) - \E_\mu (l \circ \varphi \circ \hat \sigma)\big )^2 \, d\mu(x),$$
$$V_m ^\varphi  [\hat \sigma] (l, l) : = \int _\Xx \big (l \circ \varphi\circ \hat \sigma (x) - \E_{\pb (m)} (\l \circ \varphi\circ  \hat \sigma)\big )^2\, d\pb (m) (x),$$
assuming that they   are well-defined.

%The mean square error  of a  $\varphi$-estimator $\varphi \circ \hat \sigma$  is equal to its covariance, if $\varphi \cir c\hat \sigma $ is unbiased,  i.e. $\E_\mu l \circ \varphi = l \circ \varphi d 
Let $(\Pp_\Xx, \Dd)$   be  a    2-integrable  $C^k_\tau$-diffeological statistical model.  We  call  $\varphi \circ \hat \sigma: \Xx \to V$  {\it a  regular  $\varphi$-estimator},
if for all  $l \in V'$     and  for  
for all  $\mu_0 \in \Pp_\Xx$  
\begin{equation}\label{eq:reg1}\lim_{\mu \to \mu_0} \sup \| l\circ \varphi\circ \hat \sigma  \|_{L^2 (\Xx, \mu)}  < \infty.
\end{equation}
If $\varphi\circ \hat \sigma$  is     regular, then  the function  $\varphi^l _{\hat \sigma}: (\Pp_\Xx, \Dd) \to \bR, \mu \mapsto  \E_\mu  (l \circ \varphi \circ \hat \sigma), $  is differentiable \cite[Proposition 2]{Le2020}.  
Similarly,  for a  2-integrable   parameterized statistical model $(M, \Xx, \pb)$   we call   $\varphi\circ \hat \sigma: \Xx \to V$ a   {\it regular $\varphi$-estimator}  	if for all  $l \in V'$   
for all  $m_0 \in M$ 
$$\lim_{m \to m_0} \sup \|l\circ  \varphi\circ \hat \sigma  \|_{L^2 (\Xx, \pb(m))}  < \infty.$$
If $\varphi\circ \hat \sigma$  is     regular, then  the function  $\varphi^l _{\hat \sigma}: M \to \bR, m \mapsto  \E_{\pb (m)}  (l \circ \varphi \circ \hat \sigma), $  is G\^ateaux-differentiable \cite[Lemma 5.2, p. 282]{AJLS17}.

Let $T^\g_\mu (\Pp _\Xx, \Dd)$  be the  completion  of $T_\mu (\Pp_\Xx, \Dd)$ by $\g$.  Since  $T_\xi ^\g \Pp_\Xx$  is   a Hilbert  space,  the map  
$$L_\g: T_\xi^\g \Pp_\Xx \to (T_\xi^\g \Pp_\Xx)',\, L_\g (v)(w) := \la  v, w\ra_\g,  $$   is an isomorphism. We   define  the  inverse $\g ^{-1}$ of  $\g$ on $ (T_\xi^\g \Pp_\Xx)'$ as  follows:
\begin{equation}
\label{eq:Finv}
\la   L_g  v, L_g  w \ra _{\g^{-1}}: = \la  v, w\ra_\g.
\end{equation} 
% We define $\nabla_\mu \varphi^l _{\hat \sigma}  \in T^\g_\mu (\Pp _\Xx, \Dd)$  by the formula $\la \nabla_\mu \varphi^l _{\hat \sigma} ,  V\ra  =   d_\mu  \varphi^l _{\hat \sigma}(V)$ for
Similarly, let $ {T_m^{\hat \g}} M$ be the  completion of the  quotient   space
$T_mM/\ker d_m\pb $ with  the  induced  Fisher  metric $\hat \g$. Since  the  differential  $d_m\varphi^l _{\hat \sigma}$ vanishes  on $\ker  d_m\pb$, %  there  exists   a unique   element $\nabla \varphi^l _{\hat \sigma}  \in T_m^{\hat \g} M$ such that  for any $V \in T_m M$ we have   $d_m \varphi^l _{\hat \sigma} (V) = \hat \g ({\rm pr} (V),\nabla \varphi^l _{\hat \sigma}  )$ where  ${\rm pr}$ denotes the   projection.  Note  that  
$d_m\varphi^l_{\hat \sigma}$ descends  to an  element $\hat d  _m \varphi^l_{\hat \sigma} $ in  $(T_m^{\hat \g}M)'$, which is   endowed  with  the     nondenenerate   bilinear  for $\hat  \g ^{-1}$, defined      as  in  \eqref{eq:Finv}.

The {\it diffeological  Cram\'er-Rao  inequality} asserts  that under  the  above conditions  we have \cite[Theorem 3]{Le2020}
\begin{equation}\label{eq:crd}
V_{\mu}^\varphi[\hat \sigma](l, l)-  \|d_\mu\varphi^l_{\hat \sigma}\|^2_{ \g^{-1}} \ge 0.
\end{equation}
The  {\it parametric Cram\'er-Rao} inequality asserts  that under  the  above conditions  we have \cite[Theorem  5.7]{AJLS17}, \cite{LJS2017}, \cite{JLS2017}
\begin{equation}\label{eq:crp}
V_{m}^\varphi[\hat \sigma](l, l)-  \|\hat d_m\varphi^l_{\hat \sigma}\|^2_{ \hat \g^{-1}} \ge 0.
\end{equation}
If $  M =  \Pp_\Xx $ and $\varphi: M  \to \bR^n$  is a    coordinate mapping in  a neighborhood $U(x)$ of $x \in M$  i.e. $\{x^l : = l \circ \varphi (x)\vert, l = [1, n]\}$  are  local coordinates  of $x$, then   the  diffeological   and  parametric Cram\'er-Rao  inequalities \eqref{eq:crd}, \eqref{eq:crp} become the classical Cram\'er-Rao inequality \cite[Theorem 2.2, p. 32]{AN00}:
$$	 V_m [\hat\sigma] \ge \g ^{-1}_m.$$ % \label{eq:clCR}
The proof of the diffeological and parametric   Cram\'er-Rao  inequalities is based  on 
the explicit  expression of     the dual  of $d_\mu \varphi^l_{\hat \sigma}$ in  $T^\g_\mu (\Pp_\Xx, \Dd)$, and the  dual of  $\hat d _m\varphi ^l_{\hat \sigma}$ in $T^{\hat \g}_m M$, respectively,   as the orthogonal projection of the measure $\big (l \circ \varphi \circ \hat \sigma - \E_\mu (l \circ \varphi\circ  \hat \sigma)\big ) \mu\in \Ss^{1/2}  (\Xx,\mu)_2$,  and the measure $\big (l \circ\varphi\circ \hat \sigma - \E_{\pb (m)} (l \circ \varphi\circ  \hat \sigma) \big )\pb(m)\in \Ss^{1/2} (\Xx,\pb (m))_2$, respectively,  to the  closed  subspace $ T^\g_\mu (\Pp_\Xx, \Dd) \subset  \Ss^{1/2}(\Xx, \mu)_2 $  and  the close  subspace $T^{\hat \g}_{m} M$, identified  with  the closure  of its image   via  $d_m \pb$ in $ \Ss^{1/2} (\Xx, \pb (m))_2$,   respectively.

\begin{remark}\label{rem:cr} 	
	In \cite{Janssen03}  Janssen  showed   that  $\hat \varphi ^l_\sigma: \Pp_\Xx \to \bR$   is differentiable   under  the
	weaker    assumption  that   smooth  curves  in  $(\Pp_\Xx, \Dd)$ are $L^2$-differentiable and the  regularity  of $\varphi \circ \hat \sigma$  in \eqref{eq:reg1} also  holds.  Thus his nonparametric Cram\'er-Rao inequality assumes  the weakest  condition, according  to our best  knowledge.
\end{remark}

\section{Conclusions  and final remarks}\label{sec:final}

(1) The fruitful  concept  of a  smooth   statistical  model     endowed  with the Fisher  metric and  the Chentsov  tensor   has been    investigated  and    generalized  using   different, but closely related,   formalisms   of  statistical manifolds and  natural  differentiable structures
on statistical models.  In particular,  the concept of the Fisher  metric
is naturally    related  to the   concept  of $L^2$-differentiability, the concept of a 2-integrable  parametrized  measure  model,  and  the  concept  of   statistical diffeology  $\Dd^{1,2}_{\tau_v}$.  
Since the Amari-Chentsov tensor  on  statistical  models is  a      covariant  tensor of finite   degree, it    can be      defined  and described  similarly  using      statistical diffeologies.  We refer the reader  to \cite{Amari16}  and \cite{AJLS17}  for   many  applications  of  information geometric   structures  in different  fields  of sciences.

(2)  In this  article  we omitted     a discussion   on the monotonicity  of the Fisher metric     under   Markov  kernels. The  set of all  Markov kernels   from a measurable  space $\Xx$ to  a measurable space  $\Yy$  encompasses  the   set  of  all measurable  mappings  from $\Xx$ to $\Yy$,  and as  measurable mappings     from  $\Xx$  to $\Yy$,  Markov kernels induce  (smooth) transformations  between     (smooth) statistical  models  over  $\Xx$  to (smooth) statistical models  over $\Yy$. It is well-known that  the pull back  of the Fisher  metric   under  the Markov kernel   is weaker  than   the  Fisher  metric  on the domain,   and  this property is called  the {\it   monotonicity},% see  e.g. Amari-Nagaoka  \cite[p. 31, 32]{AN00}, 
which can   be used to characterize   Fisher metrics, see  Chentsov \cite{Chentsov72}, L\^e \cite{Le2017}. The monotonicity   also justifies  the alternative  name of Fisher-Rao  metric as   the Fisher information metric  \cite{AN00}.
It  is a natural  problem to investigate      smooth families   of Markov kernels   and their  geometry. Smooth families   of measurable  mappings, e.g.   neural  networks,    and smooth families of 
Markov kernels  play  important  role  in      statistical learning  theory, in particular  in supervised learning theory.

(3) By Remark \ref{rem:lt} the Fisher  metric  generates     the   topology on statistical models  $\Pp_\Xx \subset \Ss(\Xx)$ that is     not weaker than    the topology $\tau_v$. For many  problems in statistics the weak  topology $\tau_w$ on $\Pp(\Xx)$ is more relevant. It  is an interesting problem  to  describe  statistical diffeologies  that support  metrics   generating  $\tau_w$.

(4) In  \cite[Chapter 7]{AN00}   Amari  and Nagaoka  discussed  Information  Geometry   for quantum systems, in particular  they  defined a quantum version  of the   Fisher metric   and Amari's  $\alpha$-connections    and applied them to quantum estimation theory. 
%Their  theory  has been   extended in many aspects, in particular,  using the formalism of  parameterized  measure  models in  Ciaglia, Jost  and Schwachh\"ofer \cite{CJL20}. 
It   is an  
interesting problem to develop        diffeologies   that  carry  the quantum Fisher   metric  and  the quantum  Amari's $\alpha$-connections.

\section*{Acknowledgments}
A large part of materials  in the  present  paper      was  obtained   or  digested   during    the author works  with  Nihat Ay, J\"urgen Jost and Lorenz Schwachh\"ofer on   Information Geometry, and  during the author  work with  Alexey Tuzhilin on   diffeological  Fisher  metric. She  would like    to thank them  for  the  fruitful collaboration.  She  is grateful  to Professor     Shun-ichi  Amari  for his inspiring works  and his kind
support.  She  thankfully   acknowledges   Frederic Barbaresco, Frank Nielsen, Giovani Pistone, Jun Zhang,   Patrick Iglesias-Zemmour, Enxin Wu, Jean-Pierre  Magnot,  Kaoru Ono  and   Yong-Geun Oh for stimulating   discussions on subjects   in  this  article.
A part  of this paper has been prepared during   author's visit  to the  Max-Planck-Institute for Mathematics  in Sciences in Leipzig in May 2022 and      during author's Visiting  Professorship    at the   Kyoto University from July  till October  in 2022.  The  author   would like to thank  these  institutions for their  hospitality   and excellent  working condition.   
The  research  of this paper  was supported by   GA\v CR-project 22-00091S  and  RVO: 67985840.

\section*{Disclosure  statement} The author   states that there is
no conflict of interest.

\section*{Data availability statement}  Data sharing is not applicable to this article as no data sets were generated or analyzed during the current study.


\begin{thebibliography}{9}
	
\bibitem{AJLS15} Ay N., Jost J., L\^e H.V., Schwachh\"ofer  L.,
Information geometry and sufficient statistics, {\em Probability Theory and Related Fields}   162,  {\bf 2015}, 327--364.
\bibitem{AJLS17} Ay N., Jost J., L\^e H.V., Schwachh\"ofer  L., {\em Information geometry},  Springer  Nature: Cham,  Switzerland, 2017.	 
\bibitem{AJLS18} Ay N., Jost J., L\^e H.V., Schwachh\"ofer  L., Parametrized measure models,     {\em Bernoulli} 24, {\bf 2018},   1692--1725.
\bibitem{Amari85}  Amari S., {\em Differential-Geometric Methods in Statistics}, Lecture Notes in Statistics 28, Springer-Verlag: Heidelberg, Germany,  1985.
\bibitem{AN00} Amari S., Nagaoka H. {\em Methods of Information Geometry},  Translation of Mathematical Monographs, vol. 191.  Am. Math. Soc., 2000.
\bibitem{Amari16} Amari S.,  {\em Information Geometry and Its Applications}, Applied Mathematical Sciences, vol. 194, Springer: Berlin, Germany, 2016.
\bibitem{BBM2016} Bauer, M., Bruveris, M., Michor, P., Uniqueness of the Fisher–Rao metric on the 
space  of smooth densities, {\em Bull. Lond. Math. Soc.} 48, {\bf 2016}, 499-506. 
	%\bibitem{BH11} Baez J. C.  and Hoffnung A. E., Convenient   Categories  of Smooth Spaces, {\em Trans. A.M.S.},   363, {\bf 2011},  5789-5825.
\bibitem{Bogachev2007}  Bogachev V. I.,  {\em Measure Theory I, II}, Springer, 2007.
\bibitem{Bogachev2010}  Bogachev V. I.,  {\em Differentiable Measures and the Malliavin Calculus}, AMS, Providence, Rhode  Island, 2010.
\bibitem{Bogachev2018} Bogachev V.I., {\em  Weak convergence  of measures}, Mathematical Surveys and Monographs, vol. 234, Amer. Math. Soc.:  Providence, RI, USA, 2018.
\bibitem{BS2017} Bogachev V. I. and  Smolianov O.G.,  {\em Topological  Vector Spaces and Their  Applications}, Springer, 2017.
\bibitem{Borovkov1998}  Borovkov  A. A., {\em Mathematical statistics}, Gordon and Breach Science Publishers: Amsterdam, The Nethelands, 1998.
\bibitem{Chentsov72} Chentsov N., {\em Statistical decision rules and optimal inference}, Nauka: Moscow,  Russia,   1972,   English translation  in:  Translation of Math. Monograph vol. 53, Amer. Math. Soc.:  Providence, RI, USA, 1982.
%\bibitem{CJL20} Ciaglia F. M., Jost  J., and Schwachh\"ofer  L.,  Differential geometric aspects of parametric estimation theory
%for states on finite-dimensional C$^*$-algebras. Entropy, 22(11):1332, 2020. 
%\bibitem{CJL20b} Ciaglia F. M., Jost  J., and Schwachh\"ofer  L.,  From the Jordan product to Riemannian geometries on classical and quantum states. Entropy, 22(06):637-27, 2020. 
	%\bibitem{CW2016} Christensen J.  D., Wu E., Tangent spaces and tangent bundles for diffeological spaces, {\em  Cahiers de Topologie et Geométrie Différentielle Catégoriques}, 57 {\bf 2016}, 3-50.                                                                           
	\bibitem{Fomin68} Fomin S.V.,  Differentiable measures in linear spaces. (Russian)
	{\em Uspehi Mat. Nauk} 23, {\bf 1968}, no. 1 (139), 221-222.
	\bibitem{Friedrich91} Friedrich T.,  Die Fisher-Information  und symplektische   Strukturen, {\em Math. Nachr.}  153, {\bf 1991},  273--296.
	\bibitem{GV17} Ghosal  S., van der Vaart  A. W., {\em Fundamentals of nonparametric Bayesian inference}, Cambridge Series in Statistical and Probabilistic Mathematics, vol. 44. (Cambridge University Press, Cambridge, 2017.
	\bibitem{IZ13} Iglesias-Zemmour P., {\em Diffeology}, Amer. Math. Soc.:  Providence, RI, USA, 2013.
	\bibitem{Janssen03}  Janssen  A, A nonparametric Cr\'amer-Rao inequality for estimators of statistical functional,   {\em Statistics $\&$ Probability Letters}, 64  {\bf 2003}, 347-358.
	%\bibitem{Jeffrey1946} Jeffrey H., An Invariant Form for the Prior Probability in Estimation Problems, {\em Proceedings of the Royal Society of London. Series A, Mathematical and Physical Sciences}  186, {\bf 1946},  453--461.
	\bibitem{JLS2017} Jost J. , L\^e  H. V. and  Schwachh\"ofer  L., The Cram\'er-Rao inequality  on singular  statistical models, {\em arXiv:1703.09403}.
	%\bibitem{KM1997}  Kriegl A. and  Michor P. W., {\em The  Convenient  Setting   of Global Analysis},   Amer. Math. Soc., 1997.
	%\color{black}
	%\bibitem{Lawvere1962} Lawvere W. F.,   The category of probabilistic mappings, {\bf 1962}.{\em Unpublished,	Available at } \url{hyttps://ncatlab.org/nlab/files/lawvereprobability1962.pdf}
	\bibitem{Lang99}  Lang  S. {\em Fundamentals  of Differential  Geometry}, Springer, 1999.
	\bibitem{Lauritzen75}   Lauritzen S., Statistical manifolds.  In: Differential Geometry  in Statistical Inference, Institute of Mathematical Statistics, {\em Lecture  Note-Monograph Series},  vol. 10 (1987).
	\bibitem{Le2005} L\^e  H. V., Statistical manifolds are statistical models, {\em  Journal of Geometry}  84, {\bf 2005}, 83-93.
	\bibitem{Le2017} L\^e H. V., The uniqueness of the Fisher metric as information metric, {\em Annals of Institute of Statistical Mathematics}  69, {\bf 2017}, 879-896. %,    arXiv:math/1306.1465.
	\bibitem{Le2020} L\^e  H. V., Diffeological Statistical Models, the Fisher Metric and Probabilistic Mappings, {\em Mathematics}  8(2),  {\bf 2020}, 167. % \url{https://doi.org/10.3390/math8020167}.
	\bibitem{LJS2017} L\^e  H. V.,  Jost J.  and Schwachh\"ofer  L. , The Cram\'er-Rao inequality  on singular  statistical models, {\em Proceedings of GSI 2017, LNCS 10589}, p. 552-560, Springer, 2017.
	\bibitem{LT2021}   L\^e  H. V.  and  Tuzhilin, A., Nonparametric Estimations and the Diffeological Fisher Metric. In: Barbaresco, F., Nielsen, F. (eds) Geometric Structures of Statistical Physics, Information Geometry, and Learning. SPIGL 2020. {\em Springer Proceedings in Mathematics \& Statistics}, vol 361. Springer, 2021. %\url{https://doi.org/10.1007/978-3-030-77957-3_7}.
	\bibitem{McCullagh02} McCullagh, P. What is a statistical model. {\em Ann. Stat.}  2002, 30, 1225–1310
	%bibitem{Neveu1970} Neveu J., {\em Bases Math\'ematiques du Calcul  de Probabilit\'es,deuxi\`eme edition}, Masson,  Paris, 1970.
	\bibitem{Newton2018} Newton, N. J.,  Manifolds of Differentiable Densities. {\em ESAIM: Probability and Statistics} 22 {\bf 2018}, 19-34.
	\bibitem{Pfanzagl1982}  Pfanzagl   J. {\em Contributions to  a General Asymptotic  Statistical Theory}, Lecture  Notes in Math. 13, Springer, New York, 1982.
	\bibitem{Pflug1996} Pflug  G., {\em Optimization of Stochastic Models}, Kluwer Academic, 1996.
	%\bibitem{Pflug1988} Pflug, G.,  Derivatives of probability measures - Concepts and applications to the optimization of stochastic systems. In:  {\em Lecture Notes in Control and Information Science}, {\bf 1988},  Vol. 103,  Springer  Verlag,, 252-274.
	%\bibitem{Pitcher63} Pitcher, T. S. Likelihood ratios for stochastic processes related by groups of transformations.I, II, Illinois J. Math. 7 (1963), 396-414,   Illinois J. Math. 8 (1964), 271-279.
	\bibitem{PS1985}Pistone   G. Sempi   C., {\em An infinite-dimensional geometric structure on the space of all the probability measures equivalent to a given one}, Annals of Statistics 23(1995), 1543-1561.
	\bibitem{Rao1945} Rao   C. R., Information and the accuracy attainable  in the estimation of statistical parameter.  {\em Bull. Calcutta Math. Soc.} 37 (1945), 81-89.
	%\bibitem{RSdiff}  Monthly Global Diffeology Seminar, \url{https://researchseminars.org/seminar/MGDS}.
	\bibitem{Skorohod1974} Skorohod   A.V.,  {\em Integration on Hilbert space}. Translation from Russian,
	Springer  Verlag, Berlin-New York, 1974.
	\bibitem{Souriau1980} Souriau J.-M., Groupes diff\'erentiels, {\em  Lect. Notes in Math.},{\bf  1980} vol. 836, Springer Verlag,  91--128.
	\bibitem{Strasser1985}Strasser   H. , {\em Mathematical theory of statistics}. In: De Gruyter Studies in Mathematics, Vol. 7. de Gruyter, Berlin, 1985.
	\bibitem{Vaart1991} van der Vaart   A. W., On differentiable functionals. {\em Annals of Statistics}, 19(1991).
	%\bibitem{Voiculescu93}  Voiculescu  D., The analogues of entropy and of Fisher's information measure in free probability theory. I.
	%{\em Comm. Math. Phys.} 155 (1993),  71-92.
\end{thebibliography}
\end{document}